\documentclass{llncs}
\usepackage{makeidx}  
\usepackage{float}
\usepackage{indentfirst}
\usepackage{graphicx} 
\usepackage{booktabs} 
\usepackage{caption}
\usepackage{amsmath}
\usepackage{amsfonts}
\usepackage{textcomp}
\usepackage{xcolor}

\def\R{\mathbb R}

\begin{document}
\mainmatter              
\title{On Symmetry Groups of Some Quadratic Programming Problems
\thanks
{
 The final authenticated publication is available online at
\doi{10.1007/978-3-030-49988-4}
}
}
\titlerunning{Symmetries of Quadratic Programming Problems}  
%
\author{Eremeev~A.V.\inst{1}\orcidID{0000-0001-5289-7874},
Yurkov~A.S.\inst{2}\orcidID{0000-0002-8344-0882}}
\authorrunning{Eremeev~A.V., Yurkov~A.S.} 
\institute{Sobolev Institute of Mathematics, Omsk, Russia,\\
\email{eremeev@ofim.oscsbras.ru}
 \and
Institute of Radiophysics and Physical Electronics Omsk Scientific Center SB RAS, Omsk, Russia,\\
\email{fitec@mail.ru} }

\maketitle              

\pagestyle{plain}

\begin{abstract}
Solution and analysis of mathematical programming problems may be simplified when these problems are symmetric under appropriate linear transformations. In particular, a knowledge of the symmetries may help reduce the problem dimension, cut the search space by linear cuts or obtain new local optima from the ones previously found. While the previous studies of symmetries in the mathematical programming usually dealt with permutations of coordinates of the solutions space, the present paper considers a larger group of invertible linear transformations. We study a special case of the quadratic programming problem, where the objective function and constraints are given by quadratic forms, and the sum of all matrices of quadratic forms, involved in the constraints, is a positive definite matrix. In this setting, it is sufficient to consider only orthogonal transformations of the solution space. In this group of orthogonal transformations, we describe the structure of the subgroup which gives the symmetries of the problem. Besides that, a method for finding such symmetries is outlined, and illustrated in two simple examples.

\end{abstract}

\section{Introduction}
\label{sec:intro}

Solution and analysis of mathematical programming problems may be simplified when these problems are symmetric under appropriate linear transformations. In particular, a knowledge of the symmetries may help reduce the problem dimension, cut the search space by symmetry-breaking linear cuts or obtain new local optima from the ones previously found. These methods are applicable in the case of a continuous solutions domain~\cite{CHL13,GATERMANN200495,KWM19} as well as in the integer programming~\cite{BHJ13,C99,Kolokolov2012,Margot2010,Simanchev96} and in the mixed integer programming~\cite{L10,PR19}. While most of the applications of symmetries are aimed at speeding up the exact optimization algorithms, yet in some cases the knowledge of symmetries may also be useful in designing evolutionary algorithms~\cite{Adam2004} and other heuristics.

In the present paper, we study the case of continuous solutions domain. While the previous studies of symmetries in mathematical programming usually dealt with permutations of coordinates of the solutions space~\cite{Kolokolov2012,KWM19,L10}, the present paper considers a larger group of invertible linear transformations. We study the special case of quadratically-constrained quadratic programming problem in~$\R^N$, where the objective function and the constraints are given by quadratic forms, $A, $ and $B_1,\dots,B_M$ respectively:
\begin{equation}
\label{eq:init}
\left\{
\begin{array}{l}
\displaystyle
x^TAx \to {\rm max} \, , \\
\displaystyle
x^TB_1x \le  1 \, , \\
\displaystyle
\dots \\
x^TB_Mx \le 1 \, ,
\end{array}
\right.
\end{equation}
where $x$ is an $N$-component column vector of variables, and the superscript~$T$ denotes matrix transposition. In what follows, without loss of generality we assume that $N\times N$ matrices $A, B_i, i=1,\dots,M$ are symmetric (note that any matrix can be decomposed into a sum of symmetric matrix~$S$ and  skew-symmetric matrix~$C$, and the quadratic form $x^TCx$ identically equals zero). A more substantial assumption that we will make in this paper is that $B_{\Sigma}:=\sum_{i=1}^M B_i$ is a positive definite matrix. An example of application of quadratic programming problems with such a property in radiophysics may be found e.g. in~\cite{ETY19}. 

The results of this paper may also be used for finding symmetries if some of the problem constraints have the inequality~$\le,$ some have the inequality~$\ge$ and some have the equality sign. We will consider only the inequalities~$\le$ for the notational simplicity.
The obtained results may also be applied in semidefinite relaxation methods, see e.g.~\cite{Shor1998}. Note that in~\cite{Shor1998} the well-known Maximum Cut problem (which is NP-hard) is reduced to  the problem considered here.

By a symmetry of problem (\ref{eq:init}) we mean a linear transformation
\begin{equation}
\label{eq:Lin}
x \to y=Px \, ,  
\end{equation}
defined by a non-degenerate matrix $P$ such that the problem (\ref{eq:init}), expressed in terms of the transformed space 
(i.e., through the vector columns $ y $), coincides with the original problem. That is, in terms of the vectors $ y $ our optimization problem again has the form
\begin{equation}
\label{eq:Tinit}
\left\{
\begin{array}{l}
\displaystyle
y^T A y \to {\rm max} \, , \\
\displaystyle
y^TB_1y \le  1 \, , \\
\displaystyle
\dots \\
y^TB_My \le 1 \, ,
\end{array}
\right.
\end{equation}
{\it with the same} matrix $ A $ and {\it the same set} of matrices $\{B_i: i=1,\dots,M\}$. We emphasize that, in the set of constraints, matrices $ B_i $ may be numbered arbitrarily, which, obviously, does not change the problem.
The transformations given by the matrices $P$ obviously form a group, which we denote by~$\mathcal G$. The goal of the paper is to analyse group~$\mathcal G$ and propose an algorithm for finding it.

In some cases, it may also be of interest to find the symmetry group of the set of constraints only.
Obviously, this is not much different from the search for symmetry group~$\mathcal G$ of the problem; one just needs to exclude matrix $ A $ from the consideration (i.e. formally assume that $A$ is  a zero matrix).
Furthermore, the set of symmetries of the constraints is not larger than the set of all invertible linear transformations, bijectively mapping the feasibility domain of the problem ${\mathcal D}:=\{x\in \R^N \ : \ x^T B_i x \le 1, \ i=1,\dots,M\}$ onto itself.  Therefore, the symmetry group of the set of constraints is a subgroup in the symmetry group of invertible linear transformations of~${\mathcal D}$.

The structure of the paper is the following. In Section~\ref{sec:symmetry}, it is shown that the group of linear symmetries of the problem is a subgroup of orthogonal transformations. Also, the structure of the group of symmetries and the corresponding Lie algebra are discussed. In Section~\ref{sec:algo}, a general algorithm for finding the symmetries is proposed, and in Section~\ref{sec:examp} it is illustrated in two simple examples. A discussion of the results and the conclusion are in Sections~\ref{sec:disc} and~\ref{sec:conc}. 
Appendix contains a proof of a ``folklore'' fact from matrix analysis used in Section~\ref{sec:algo}.

\section{Structure of the Symmetry Group}
\label{sec:symmetry}

Invariance of the problem under transformation~$P$ implies that
\begin{equation}
\label{eq:geninvar}
P^TAP=A \,  , \qquad P^TB_i P = \sum_{j=1}^M L_{ij}B_j , \  j=1,\dots,M,
\end{equation}
where $ L_{ij} $ are the elements of a permutation matrix, i.e. matrix ${L=(L_{ij})}$ has a single ``1'' in each column and in each row, other elements of~$L$ are zeros. 

If (\ref{eq:geninvar}) holds, then the invariance condition of the matrix~$B_{\Sigma}$ is satisfied:
\begin{equation}
\label{eq:invBS}
P^TB_{\Sigma}P = B_{\Sigma}.
\end{equation}
Naturally, the converse is not true in the general case, but at least we can say that the desired group~$\mathcal G$ is a 
subgroup of the invariance group of matrix $B_{\Sigma}$.
This matrix may be represented as a congruent transformation of a diagonal matrix:
\begin{equation}
\label{eq:BSSTDS}
B_{\Sigma}=S^TDS \, ,
\end{equation}
where $D$ is a diagonal matrix, which can have only ``0'', ``1'', or ``-1'' on its main diagonal. Essentially, we are talking about reducing the quadratic form corresponding to matrix $B_{\Sigma} $ to its canonical form. So  matrix $S$ can be constructively obtained, for example, by the finite Lagrange method (\cite{Lancaster}, Ch.~5). 

Now, if we restrict ourselves to the special case where matrix $B_{\Sigma} $ is positive definite (it occurs, for example, in the radiophysical problem of optimizing the excitation of antenna arrays~\cite{ETY19}), then $ D $ will be the unit matrix and it may be omitted in~(\ref{eq:BSSTDS}). Condition~(\ref{eq:invBS}) then turns into
\begin{equation}
P^TS^TSP =  S^TS \, 
\end{equation}
or
\begin{equation}
(SPS^{-1})^T  (SPS^{-1}) =  E \, ,
\end{equation}
where $E$ is a unit matrix. This means that matrix
\begin{equation}
\label{eq:Qdef}
Q:=SPS^{-1}
\end{equation}
is in the group of orthogonal transformations~$O(N)$ (see e.g.~\cite{Zhelobenko}). So we proved
\begin{proposition} \label{prop:subgroup_of_ON}
If $B_{\Sigma}$ is positive definite then group~$\mathcal G$ is isomorphic to some subgroup of~$O(N)$ and this isomorphism
is given by Equation~(\ref{eq:Qdef}).
\end{proposition}

Since $P=S^{-1}QS$ by~(\ref{eq:Qdef}), so application of~(\ref{eq:geninvar}) gives
\begin{equation}
\begin{array}{l}
\displaystyle
(S^{-1}QS)^TA(S^{-1}QS)=A \,  , \\ \\
\displaystyle
(S^{-1}QS)^TB_i (S^{-1}QS)= \sum_{j=1}^N L_{ij}B_j \, , \ i=1,\dots,M,
\end{array}
\end{equation}
and after a simple transformation we have
\begin{equation}
\label{eq:invarwithQ}
Q^T \tilde{A} Q=\tilde{A} \,  , \ \ \ \
Q^T \tilde{B_i} Q= \sum_{i=1}^N L_{ij}\tilde{B_j} \, , \ i=1,\dots,M,
\end{equation}
where
\begin{equation}
\tilde{A}=\left(S^{-1}\right)^T A S^{-1} \, , \qquad
\tilde{B_i}=\left(S^{-1}\right)^T B_i S^{-1} \, , \ i=1,\dots,M.
\end{equation}
So using isomorphism (\ref{eq:Qdef}) we can substitute equations (\ref{eq:geninvar}) by the similar equations  (\ref{eq:invarwithQ}), but with the matrix substitution
\begin{equation}
\label{eq:totilde}
A \to \tilde{A} \, , \qquad B_i \to \tilde{B_i}\, , \ i=1,\dots,M. 
\end{equation}
and substituting $P$ by the orthogonal matrix~$Q$. These equations are significantly simpler, since in this case condition~(\ref{eq:invarwithQ}) may be formulated linearly in~$Q$:
\begin{equation}
\label{eq:comutatinvarwithQ}
\tilde{A} Q=Q\tilde{A} \,  , \ \ \ \
\tilde{B_i} Q= Q\sum_{j=1}^M L_{ij}\tilde{B_j} \, , \  i=1,\dots,M.
\end{equation}
If one finds all suitable orthogonal mappings $ Q$, then it will be easy to restore the corresponding matrices~$P $. 
Assuming all this, we omit the tildes above matrices~$A$ and $B_i$ further in order to simplify the notation.

It is well-known that the orthogonal group $O(N) $ consists of two connected components, for one of them the determinant of the matrix equals~1, for the other it equals~-1 (see e.g.~\cite{Zhelobenko}). The first component is a subgroup of $O(N)$, denoted by $SO(N)$ and also called the {\em rotation group,} due to the fact that in dimensions 2 and 3, its elements are the usual rotations around a point or a line, respectively. The second component does not constitute a subgroup of $O(N)$, since it does not contain the identity element. Matrices from the second component can be represented, for example, in the following form: $ {\rm diag} \{- 1,1 \dots 1 \} Q $, where $ Q \in SO(N) $, so between these components there is a one-to-one correspondence (which is not an isomorphism in the group-theoretical sense, since it does not preserve the group operations). The  required matrices $Q$ can belong to both the first component and the second.

The standard facts of topological groups theory (see e.g.~\cite{Zhelobenko},~Ch.~1) imply the following properties of symmetry group~$\mathcal G,$ endowed with the standard topology of $\R^{N^2},$ applicable to the space of $(N\times N)$-matrices. As any topological group, $\mathcal G$  consists of connected components (in the topological sense), only one of which, hereafter denoted as~$\mathcal G_1$,  contains the identity element. This~$\mathcal G_1$   is invariant subgroup of ~$\mathcal G$, see Theorem~1~\cite{Zhelobenko},  and   called {\em the continuous subgroup of symmetries} in what follows.  The remaining connected components (not being subgroups) can be considered as products of the elements of the group outside~$\mathcal G_1$ and the elements of~$\mathcal G_1$ i.e. the cosets of~$\mathcal G_1$. These cosets make up a discrete group. Given that~$\mathcal G_1$ is an invariant subgroup, multiplication of the cosets of this discrete group is determined naturally, and the discrete group is  a factor group $\mathcal G / \mathcal G_1$. These cosets can be identified by indicating one (any) representative of a coset.

Naturally, degenerate cases are possible. First, when ~$\mathcal G_1$  degenerates into the identity element, the entire symmetry group~$\mathcal G$ is a purely discrete group (each coset consists of one element). Secondly, there may be no other elements of discrete symmetry but  only the continuous subgroup of symmetries~$\mathcal G_1$  . And finally, the entire symmetry group~$\mathcal G$ may consist of only the identity element.

\section{Finding the Symmetry Group}
\label{sec:algo}

Due to the observations from Section~\ref{sec:symmetry}, the search for all appropriate symmetry transformations $Q $ may be divided into two parts: the search in the first component of $ O(N) $ (i.e., in subgroup~$SO(N) $) and the search in the second component where the determinant of orthogonal matrices equals~-1. 
Initially  we  restrict ourselves  to the first subset.   A   generalization to the whole group $O(N) $  will  be done by analogous consideration of the second subset while searching for discrete symmetries. The only difference will be that in the second case, it will be necessary to replace $ Q \to {\rm diag} \{- 1,1, \dots ,1 \} Q $.

\subsection{Continuous Subgroup of Symmetry}
First, we consider the continuous subgroup of symmetry~$\mathcal G_1$. Nontrivial permutations of matrices $ B_i$ can not result from transformations which belong to ~$\mathcal G_1$ , since it is impossible to continuously move from the identical transformation (which implies that matrices $B_i$ are not permuted) to any transfomration~$Q$ yielding a non-trivial permutation of matrices $B_i$. Note that any such~$Q$ has a neighborhood of transformations which do not yield the trivial permutation of the matrices~$B_i$. So the invariance conditions must hold:
\begin{equation} \label{eqn:matrix_transform}
Q^T AQ = A  \, , \qquad  Q^TB_iQ = B_i\, , \ i=1,\dots,M.
\end{equation}  
For orthogonal transformations $Q $, this is equivalent to commutativity:
\begin{equation}
\label{eq:commutat}
A Q = QA\, , \qquad B_i Q = QB_i\, , \ i=1,\dots,M.
\end{equation}
Proposition~(\ref{prop:skew_exp}) is a ``folklore'' fact of matrix analysis (the proof is in appendix):

\begin{proposition} \label{prop:skew_exp}
Any matrix  $Q \in SO(N) $ can be represented as a matrix exponential function of a skew-symmetric matrix. The converse is also true: the exponential function of any skew-symmetric matrix is an orthogonal matrix.
\end{proposition}

So with some skew-symmetric matrix~$X$ we have $Q=e^X$. The set of skew-symmetric matrices $X$ make up the Lie algebra corresponding to  this    Lie group~\cite{Zhelobenko}. (The Lie algebra corresponding to $ SO(N) $ is usually denoted by $ so(N) $.) Any Lie algebra is also a linear space, any of its elements can be expressed by means of basis elements, called generators. Thus, any element of the Lie algebra can be represented as:
\begin{equation}
X= \sum_n a_n G_n \, ,
\end{equation}
where $ a_n $ are real numbers, $ G_n $ are the generators. The space of skew-symmetric matrices has a dimension $ N (N-1) / 2 $, and there will be as many coefficients $ a_n $ and as many generators. As generators, one can choose matrices containing one unit element above the main diagonal (the rest are zeros), then the skew-symmetry uniquely determines the remaining matrix elements of these generators. So, any element $ Q $ of $ SO (N) $ can be represented as:
\begin{equation}
\label{eq:sunexp}
Q=e^{\sum\limits_n a_n G_n} \, .
\end{equation}

Since the desired continuous subgroup of symmetry ~$\mathcal G_1$   is a subgroup of $ SO (N) $, so representation~(\ref{eq:sunexp}) is also valid for it, but, generally speaking, the parameters $ a_n $ are not independent now. Thus, the search for this subgroup essentially reduces to finding the restrictions on parameters $ a_n $.

It is quite obvious that in order for commutativity conditions~(\ref{eq:commutat}) to be satisfied, it suffices that the following conditions hold true:
\begin{equation}
\label{eq:commutat2}
\left\{
\begin{array}{l}
\displaystyle
B_i \left(\sum\limits_na_nG_n\right) = 
\left(\sum\limits_na_nG_n\right)B_i \, , \\ \\
\displaystyle
A \left(\sum\limits_na_nG_n\right) = \left(\sum\limits_na_nG_n\right)A \, .
\end{array}
\right.
\end{equation}
It means that if  matrix $ X $  commutes with all matrices $ B_i $ and with matrix $ A $, then $X$ lies in Lie algebra of~$\mathcal G_1$. 
Indeed, expanding the exponential function in a power series, we see that if the matrices $ A $ and $ B_i $ commute with the argument of this function, then they commute with the exponential function itself. 
Note that condition~(\ref{eq:commutat2}), generally speaking, is not necessary to fulfill (\ref{eq:commutat}). 
%
%
%
%
However, the continuous subgroup of symmetry, as a connected Lie group, is completely determined by its Lie algebra, so it is completely determined by the restrictive relations for elements of the Lie algebra\footnote{  For abstract groups, such a unique connection exists only in the case of simply connected groups; otherwise, an abstract exponent cannot be uniquely determined. But in our case of a matrix group, the matrix exponent is uniquely determined. }. Thus, in search for the continuous subgroup of symmetry, (\ref{eq:commutat}) may be replaced with (\ref{eq:commutat2}).

Equations (\ref{eq:commutat2}) are a system of linear algebraic equations that determine parameters $a_n $. This system is homogeneous, so it has a continuum of nonzero solutions. Note that there is always a trivial zero solution to the system of equations (\ref{eq:commutat2}) corresponding to an identity matrix $Q$. Some of parameters $ a_n $ remain ``free'' (these will be the parameters of the desired subgroup), and the rest of~$a_n $ may be linearly expressed through the ``free'' ones. The solution to this system of equations~(\ref{eq:commutat2}) can be obtained constructively by the Gauss method. 

The condition of problem invariance under transformation~$Q$ turnes into
\begin{equation}
\label{eq:subG1}
Q=e^{\sum_n a_n \hat{G}_n} \, ,
\end{equation}
where the sum goes over the ``free'' parameters $ a_n $, and the new generators denoted by~$\hat{G} _n$ are linear combinations of the former generators~$G_n$. The set of all $Q$ matrices satisfying~(\ref{eq:subG1}) is parameterized by a finite set of real parameters $ a_n $. Note, however, that this set of matrices is not necessarily isomorphic to a Euclidean space, since more than one set of parameters $ a_n $ can correspond to the same~$Q$.  

Let us show that the set of matrices defined by formula (\ref{eq:subG1}) is a group. To this end, it is sufficient to prove that this matrix set    $ \hat{\mathcal A} = \{ \hat{X}, \hat{X}= \sum_n a_n \hat{G}_n \}$   is a Lie algebra. For a matrix algebra to be a Lie algebra, it is necessary and sufficient to be closed relative to the calculation of the commutator,    i.e. $\hat{\mathcal A}$ is Lie algebra if and only if  for any  $ \hat{X}_i, \hat{X}_j \in \hat{\mathcal A}$ a commutator 
\begin{equation}
[\hat{X_i},\hat{X_j}] = \hat{X_i}\hat{X_j} - \hat{X_j}\hat{X_i} \, ,
\end{equation}
is also an element of   $\hat{\mathcal A}$  . This is easily verified in our case.  Indeed, since all $\hat{X} $ lie in $ so(N) $, their commutators also lie in $ so(N) $. Therefore, for them to lie not only in $so(N) $, but also in  $\hat{ \mathcal A}$   , that is, for this algebra to be a Lie algebra, it is sufficient that these commutators satisfy the same restrictive conditions that distinguish set  $ \hat{\mathcal A} $    from $ so(N) $. The restrictive conditions~(\ref{eq:commutat2}) mean that all $ \hat{X} $ commute with all matrices $ B_i$ and with matrix $ A $. But then all the products of such $ \hat{X} $ also commute with all matrices $ B_i $ and with matrix $ A $. And then the commutator $ [\hat{X_j}, \hat{X_j}] $, which is a product difference, satisfies the same restrictive conditions.  Thus, the set of matrices   $ \hat{\mathcal A} $   is a Lie algebra, and therefore the set of matrices defined by formula~(\ref{eq:subG1}) is a Lie group. 

Now let us prove that the set of matrices defined by formula~(\ref{eq:subG1}) is {\em the whole} continuous subgroup of symmetries~$\mathcal G_1$. We will show that a converse leads to a contradiction. Indeed, the converce assumption implies that in the algebra of the continuous group of symmetry there is at least one more generator $G_{\rm extra}$ (with its own coefficient, let it be~$b$), linear independent from generators~$\hat{G} _n$.  But then there is a one-dimensional subgroup of~$\mathcal G_1$ produced by the element $Q=e^{b G_{\rm extra}}$. If we substitute this~$Q$ into the invariance condition~(\ref{eq:commutat}), differentiate with respect to $b$ and set $b = 0$,   then it   turns out that $G_{extra}$ satisfies exactly the same condition, which distinguishes the set of matrices   $\hat{\mathcal A}$   from the entire Lie algebra of  group $SO(N).$ So this additional generator lies in the linear hull of the generators~$\hat{G} _n$. Which is a contradiction.
So we have proved the following

\begin{theorem}\label{thm:find_continuous}
The continuous subgroup of symmetries~$\mathcal G_1$ consists of orthogonal transformations with matrices expressed by the matrix exponential function $e^{\sum_n a_n \hat{G}_n},$ where $a_n$ are any real-valued parameters, and all $\hat{G}_n$ make up a basis of the space of solutions to the system of linear equations~(\ref{eq:commutat2}) in the linear space of the $(N\times N)$ skew-symmetric matrices.
\end{theorem}

\subsection{Discrete Group of Symmetry}

In the case of discrete symmetry, nontrivial permutations of matrices $ B_i $ are possible. Therefore, the condition (\ref{eq:commutat}) is replaced by the following:
\begin{equation}
\label{eq:transpcommutat}
AQ=QA \,  , \qquad B_i Q = Q\sum_{j=1}^N L_{ij}B_j \, , \ i=1,\dots,M.
\end{equation}
There are $M!$ permutation matrices $L$ and they can be enumerated for small problems. Then we can assume that in (\ref{eq:transpcommutat}) $ L_ {nm} $ are known. (Note that if we generalize  Problem~(\ref{eq:init}) so that some of the constraints have inequalities~$\le,$ some have inequalities~$\ge,$ and some have equalities, then the permutations in each of these three subgroups should be considered.)
Furthermore, iterating over all possible matrices $ L $, one can solve equations (\ref{eq:transpcommutat}) with respect to $Q $. But it must be taken into account that matrix $Q $ lies in $ SO(N) $, otherwise equation (\ref{eq:transpcommutat}) is not valid.  To this end, one can represent $ Q $ as a matrix exponential function~(\ref{eq:sunexp}) and solve the equation for $N(N-1)/2$ parameters $ a_n $ as variables. The same should be done with matrix $Q{\rm diag}\{-1,1,\dots,1\}$. The resulting equations will involve exopnential functions, so for their solution in each particular case, it is necessary to develop a special numerical method.
Alternatively, one can solve equations~(\ref{eq:transpcommutat}) for matrix~$Q$ as a variable, conditioned that $QQ^T=E.$

\section{Illustrative Examples}
\label{sec:examp}

\subsection{Example with Trivial Continuous Subgroup of Symmetries} \label{subsec:ex1}

\begin{figure}[h]
    \centering
            \center{\includegraphics[width=0.5\textwidth]{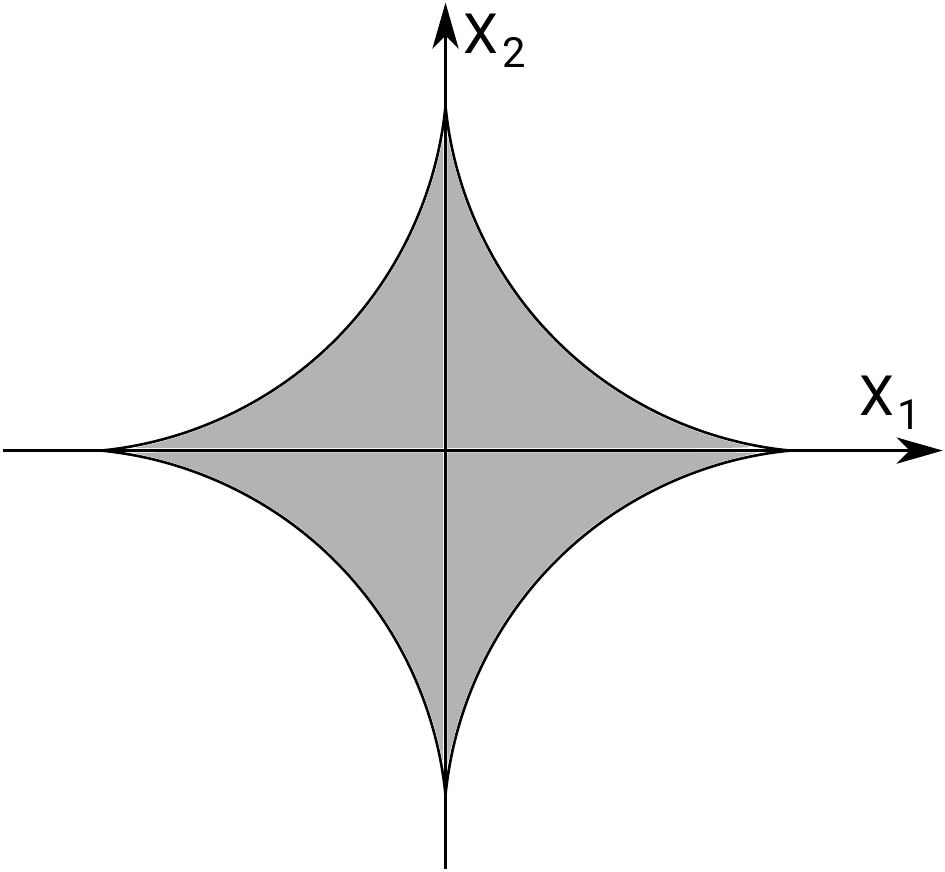} }
\caption{Feasibility domain defined by matrices $B_1$ and $B_2$ in Subsection~\ref{subsec:ex1}.}
\label{Region}
    \end{figure}

Let us apply the obtained results to a quadratic programming problem with $N=M=2,$ defined by the following matrices (see Fig.~\ref{Region})
\begin{equation}
A=\left(
\begin{array}{cc}
1.0 & 0.0 \\
0.0 & 0.8
\end{array}
\right) \, ,
\end{equation}
\begin{equation}
B_1=\left(
\begin{array}{cc}
0.5 & 2.0 \\
2.0 & 0.5
\end{array}
\right) \, ,
\ 
B_2=\left(
\begin{array}{cc}
0.5  & -2.0 \\
-2.0 & 0.5
\end{array}
\right) \, .
\end{equation}
In this example, $ B_{\Sigma}$ is the identity matrix, and so  $S=E$,  therefore  transformation~(\ref{eq:totilde}) is not necessary.
The feasibility area corresponding to matrices $ B_1 $ and $ B_2 $ is shown in Figure~\ref{Region}. Its rotational symmetry properties (as well as the symmetry properties of the problem which involves matrix $ A $) are obvious from geometric considerations: the symmetry group of the domain~$\mathcal D$ consists of the identical transformation (the identity matrix), rotations of 90~degrees, 180~degrees and 270~degrees (the latter is also the inverse element to the rotation of 90~degrees). In total, there are four elements of the group.

For the symmetry group of the problem, 90 and 270 degrees rotations disappear, the two other elements of the group remain. It is also clear that there will be four local optima, two of which are global.

Let us now verify that the results described above give the same result. Firstly, in this two-dimensional case there is only one generator:

\begin{equation}
G=\left(
\begin{array}{cc}
0 & 1 \\
-1 & 0
\end{array}
\right) \, .
\end{equation}

Accordingly, there is only one coefficient~$a $. The generator $ G $ does not commute with any of the matrices written above. Therefore, the system of equations~(\ref{eq:commutat2}) has only one zero solution corresponding to an identity matrix~$E $. The continuous subgroup of symmetry in this example degenerates into a trivial subgroup of one identity element.

To find a discrete symmetry by direct calculations, we note that
\begin{equation}
e^{aG}=
\left(
\begin{array}{cc}
\cos a  & \sin a \\
-\sin a  &  \cos a
\end{array}
\right) \, ,
\end{equation}
\begin{equation}
\begin{array}{l}
\displaystyle
B_1e^{aG}=0.5 \cos a \, E - 2\sin a \, D + 0.5 \sin a \, G
+ 2\cos a \, H \, , \\ \\
\displaystyle
e^{aG}B_1=0.5\cos a \, E + 2\sin a \, D + 0.5\sin a \, G
+ 2\cos a \, H \, , \\ \\
\displaystyle
B_2e^{aG}=0.5\cos a \, E + 2\sin a \, D + 0.5\sin a \, G
- 2\cos a \, H \, , \\ \\
\displaystyle
e^{aG}B_2=0.5 \cos a \, E - 2\sin a \, D + 0.5 \sin a \, G
- 2\cos a \, H \, , \\ \\
\displaystyle
Ae^{aG}= 0.9 \cos a \, E + 0.1 \cos a \, D + 0.9 \sin a \, G 
+ 0.1 \sin a \, H \, , \\ \\
\displaystyle
e^{aG}A= 0.9 \cos a \, E + 0.1 \cos a \, D + 0.9 \sin a \, G 
- 0.1 \sin a \, H \, .
\end{array}
\end{equation}
where
\begin{equation}
\begin{array}{l}
\displaystyle
D= \left(
\begin{array}{cc}
1  &  0 \\
0  & -1
\end{array} 
\right) \, , \\ \\ 
\displaystyle
H= \left(
\begin{array}{cc}
0  &  1 \\
1  &  0
\end{array} 
\right) \, . 
\end{array}
\end{equation}
Substituting this all into the equations from Section~\ref{sec:algo}, we obtain the following. When considering the symmetry of~$\cal D$ without permutations of matrices $B_i $, we obtain the equation $ \sin a = 0 $, and with permutations, the equation $ \cos a = 0 $. The first one corresponds to the identical transformation and a rotation of 180 degrees ($ a = 0, \pi $). The second one corresponds to rotations of 90 and 270 degrees ($ a = \pi / 2, 3 \pi / 2 $). Thus, a formal application of the above formulas agrees with the geometric considerations.

If we additionally require the symmetry of the objective function, then in both cases (with the permutation and without it) the second equation $ \sin a = 0 $ will appear, excluding rotations of 90 and 270 degrees. Finally, to obtain {\em all} symmetries of the problem, one has to solve equations~(\ref{eq:transpcommutat}) for the matrix ${\rm diag}\{-1,1\}e^{aG}$ and join the resulting symmetries with the  rotations found before.

\subsection{An Example with Non-Trivial Continuous Subgroup of Symmetries} \label{subsec:ex2}

As a second example, now with a continuous symmetry, we can take a problem with $N=3, M=2$ defined by the following matrices
\begin{equation}
\begin{array}{l}
\displaystyle
A={\rm diag}\{1,1,1\} \, , \\ \\
\displaystyle
B_1={\rm diag}\{2,2,0\} \, , \\ \\
\displaystyle
B_2={\rm diag}\{-1,-1,1\} \, . \\ \\
\end{array}
\end{equation}
In this example, the objective function is obviously invariant under any transformations from $SO(3) $, so the symmetry of the problem coincides with the symmetry of~$\mathcal D$. Again,  transformation~(\ref{eq:totilde}) is not necessary here, since $B_{\Sigma}$ is the identity matrix.

In this example, we will choose the generators in the following form:
\begin{equation}
G_1=\left(
\begin{array}{ccc}
0 &  0 & 0 \\
0 & 0 & -1 \\
0 & 1  & 0
\end{array} \right) ,  \,
G_2=\left(
\begin{array}{ccc}
0 &  0 & 1 \\
0 & 0 &  0 \\
-1& 0  & 0
\end{array} \right)  ,  \,
G_3=\left(
\begin{array}{ccc}
0 & -1 & 0 \\
1 & 0 &  0 \\
0 & 0  & 0
\end{array} \right)  .
\end{equation}
Substituting this into (\ref{eq:commutat2}) we see that $ a_1 = a_2 = 0 $, and the arbitrary parameter is $ a_3 $. Thus, the continuous symmetry subgroup is described by the following one-parameter matrix family:
\begin{equation}
\label{eq:comp1}
e^{a_3G_3} = \left(
\begin{array}{ccc}
\cos a_3 & -\sin a_3 & 0 \\
\sin a_3  & \cos a_3  & 0 \\
0 & 0 & 1
\end{array}
\right) \, .
\end{equation}
To find the discrete symmetry in this particular case, it is more convenient to represent matrix $ Q $ not in the exponential form~(\ref{eq:sunexp}) but rather through the Euler parameters $ \alpha $, $ \beta $ and $ \gamma $, as a product of three exponential functions:
\begin{equation}
\label{eq:3exp}
Q=e^{\alpha G_3} e^{\beta G_1} e^{\gamma G_3}  .
\end{equation}
Now we substitute~(\ref{eq:3exp}) into equation~(\ref{eq:comutatinvarwithQ}), which may be written as
\begin{equation}
e^{-\gamma G_3} e^{-\beta G_1}e^{-\alpha G_3} 
B_i e^{\alpha G_3} e^{\beta G_1} e^{\gamma G_3} 
= \sum_{j=1}^M L_{ij}B_j \, , i=1,\dots,M.
\end{equation}
Note that $\exp(\alpha G_3)$ commutes with both matrices $B_1, B_2$, and therefore the left factor cancels out. The last factor also cancels out after multiplying the equations on the left and on the right side by the similar exponential functions. So the defining equation~(\ref{eq:comutatinvarwithQ}) reduces to 
\begin{equation}
e^{-\beta G_1} B_i e^{\beta G_1} = \sum_{j=1}^M L_{ij}B_m \, , \ i=1,\dots,M.
\end{equation}
We have two options for permutations: one trivial and one non-trivial. Accordingly, two options are obtained. The first:
\begin{equation}
\left\{
\begin{array}{l}
\displaystyle
e^{-\beta G_1} B_1 e^{\beta G_1} = B_1 \, , \\ \\
\displaystyle
e^{-\beta G_1}B_2 e^{\beta G_1} = B_2 \, .
\end{array}
\right.
\end{equation}
the second:
\begin{equation}
\left\{
\begin{array}{l}
\displaystyle
e^{-\beta G_1} B_1  e^{\beta G_1} = B_2 \, , \\ \\
\displaystyle
e^{-\beta G_1} B_2 e^{\beta G_1} = B_1 \, .
\end{array}
\right.
\end{equation}
We note that due to the equality $ B_2 = E-B_1 $, in both cases the second equation can be reduced to the first one and vice versa. So from two equations it is enough to solve only one. By direct calculations we obtain the following:
\begin{equation}
e^{\beta G_1}=\left(
\begin{array}{ccc}
1     & 0 & 0 \\
0     & c & - s \\
0     & s & c
\end{array}
\right) \, ,
\end{equation} 
where for simplicity of notation we denote
$s:=\sin \beta, \,  c=\cos \beta \, $.
Further direct calculations give
\begin{equation}
e^{-\beta G_1} B_1 e^{-\beta G_1} = 2\left(
\begin{array}{ccc}
   1   &     0           &   0  \\
   0    &    c^2       &   -cs   \\  
   0    &    -cs         &   s^2
\end{array}
\right) \, .
\end{equation}
In the case of the trivial permutation, this reduces to a system of equations that has two obvious solutions: $ c = \pm 1 $, $ s = 0 $.
This results in two options for matrix $ Q $:
\begin{equation}
\begin{array}{l}
\displaystyle
e^{\alpha G_3}\, {\rm diag}\{1,1,1 \} \, e^{\gamma G_3} \, , \\ \\
\displaystyle
e^{\alpha G_3}\,{\rm diag}\{1,-1,-1 \} \, e^{\gamma G_3} \, .
\end{array}
\end{equation}
Obviously, the first matrix belongs to a continuous subgroup of symmetry, it does not need to be taken into account, since such matrices are already taken into account above. The second matrix, however, does not belong to the continuous subgroup\footnote{This is because $ Q_{33} $ is~$-1$, rather than~1 as in the continuous subgroup.}. As a representative of this component, we can take the above expression, written for $\alpha = \gamma = 0 $, i.e. just $ {\rm diag} \{1, -1, -1 \} $.

In the second case, where the permutation of matrices $B_i$ is non-trivial, the system of equations obviously has no solutions.

Thus, the  subgroup of orthogonal symmetries with determinant~1 in this example consists of two connected components. The first one is described by the matrix 
family~(\ref{eq:comp1}), parametrized by one real parameter (rotation angle). The second one is described by the same matrices, but multiplied by ${\rm diag} \{1, -1, -1 \} $.

To obtain the whole group~$\mathcal G$, one has to solve the equations from Section~\ref{sec:algo} for the matrix ${\rm diag}\{-1,1,1\}Q$ and join the resulting symmetries $Q$ to the subgroup of orthogonal symmetries with determinant~1 which we found above.

\section{Discussion}
\label{sec:disc}

As a ``brute force'' approach to finding all symmetries of the problem, one can formulate a non-linear optimization problem
 in $\R^{N^2}$: 
$$
 \min\left( ||A Q - QA|| + \sum_{i=1}^M || B_i Q - Q C_i(L) || \ : \ QQ^T=E\right)
$$
where $Q$ is a matrix of variables, the matrices $C_i(L)$ are defined by~$L$ as $C_i(L):=\sum_{j=1}^M L_{ij}{B_j}$, and $||\cdot||$ denotes any matrix norm. A set of optimal solutions (with zero objective value) gives the set of orthogonal symmetry transformations. The union of $M!$ such sets, taken over all permutation matrices~$L$, makes up the whole group~$\mathcal G$. In the case of trivial continuous subgroup of symmetry, each of the $M!$ problems has a discrete set of optimal solutions, which, in principle, may be found e.g. by a multi-start of a gradient descent method.

There are other options to find group~$\mathcal G$ using non-linear programming. For example, one can similarly formulate
a minimization problem with respect to the elements of matrix~$P$. Moreover,
there is no need to impose the condition $\det(P) \ne 0$, since it follows from~(\ref{eq:invBS}) that the square of this determinant is equal to~1.
Analysis of the properties and methods of solution of these non-linear optimization problems is beyond the scope of the paper.

In applications of quadratic programming, it is not necessary to find all symmetries of a problem to improve performance of solution algorithms, such as the branch and cut method. If some valid cuts are known already for the problem instance, then each linear symmetry of a problem may be used to double the set of valid cuts. Even if there were no cuts known before, then any symmetry $P\in \mathcal G$ which maps a hemi-space $\{x : a^Tx \ge 0\}$ into the  hemi-space $\{x : a^T x \le 0\}$ with some $a\in \R^N$ then the constraint~$a^Tx \ge 0$ may be added to the set of problem constraints as a valid cut. 

If $\mathcal G$ has a non-trivial continuous subgroup so large that any element of~$\mathcal D$ may be mapped onto some hyper-plane in~$\R^N$ by a corresponding $P \in \mathcal G$, then the problem dimension may be decreased by one, see e.g. the problem from Subsection~\ref{subsec:ex2}, where any vector may be rotated by mapping~(\ref{eq:comp1}) with an appropriate angle~$a_3$ into the subspace~$\{x: x_1=0\}$.  In this respect, it would be appropriate to study the following hypothesis: {\em Problem~(\ref{eq:init}) may be reduced to a problem of the same form in solutions space $\R^{N-K},$ where $K$ is the size of the basis mentioned in Theorem~\ref{thm:find_continuous}.}

In local search, the problem symmetries may be used to identify equivalence classes of local optima (consisting of local optima, identical up to a symmetry transformation) since obviously, local optima are mapped to local optima under invertible linear symmetries of the problem. In the multi-start procedure, a smaller number of visited equivalence classes, compared to the number of visited local optima, should tighten estimates of the total  local optima number~\cite{GK,RE04}.

\section{Conclusions}
\label{sec:conc}
The results obtained in this paper further extend the applicability of the approach to improving algorithms performance in the mathematical programming, employing symmetries of the problem. The authors are not aware of other works on problem symmetries, based on the theory of Lie groups and Lie algebras. It is expected that the proposed approach may be extended to other types of problems in the mathematical programming. In particular, it would be interesting to try extending the analysis to the general case of problem~(\ref{eq:init}) without the assumption of positive-definiteness of the sum of matrices of quadratic forms. It is challanging in this case that instead of the group of orthogonal transformations~$O(n)$ one would have to consider the more general pseudo-orthogonal group~$O(p,q)$. Technical development of the outlined method for finding problem symmetries is also a subject of further research.

\section*{Acknowledgment}

The authors thank V.M. Gichev for helpful comments on the preliminary version of the manuscript.
The work on Sections~\ref{sec:symmetry} and~\ref{sec:algo} was funded in accordance with
the state task of the Omsk Scientific Center SB RAS (project
AAAA-A19-119052890058-2).
The rest of the work was funded by the program of fundamental scientific research of the SB RAS, I.5.1., project 0314-2019-0019.

\bibliographystyle{splncs03}
\bibliography{MOTOR2020}

\section*{Appendix}

This appendix contains the proof of Proposition~\ref{prop:skew_exp} and will be removed from the final version of the paper in the case of acceptance. To prove Proposition~\ref{prop:skew_exp}, we will use the following statement. 

\begin{proposition} \hspace{-1em} \label{prop:commute}
\footnote{See, e.g.~\S~2.7 in Onishchik, A.L., Vinberg, E.B.: Lie Groups and Algebraic Groups. Springer (1990)}
If two matrices $ M_1 $ and $ M_2 $ commute, i.e. $M_1 M_2=M_2 M_1,$  then $ e^{M_1} e^{M_2} = e^{M_1 + M_2}$.
\end{proposition} 

\noindent
{\bf Proof of Proposition~\ref{prop:skew_exp}.}
The orthogonal matrix is normal. 
Therefore\footnote{According to Theorem~2.5.8 from Horn, R.A., Johnson, C.R.: Matrix Analysis. Cambridge University Press (2012)}, 
the matrix $ M $ is orthogonally equivalent to a block-diagonal matrix $ A $, i.e.:
\begin{equation}
Q=W^T A W \, .
\end{equation}
$ W $ is an orthogonal matrix, the blocks of the matrix $ A $ have a size of no more than $ 2 \times 2 $. In particular, the blocks of size $ 2  \times 2 $ have the form:
\begin{equation}
\label{eq:block}
\left(
\begin{array}{cc}
a  & b \\
-b & a
\end{array}
\right) \, .
\end{equation}
The orthogonality of $ Q $ implies the orthogonality of such blocks, i.e.
\begin{equation}
\left(
\begin{array}{cc}
a  & b \\
-b & a
\end{array}
\right)^T
\left(
\begin{array}{cc}
a  & b \\
-b & a
\end{array}
\right) = 
\left(
\begin{array}{cc}
a^2+b^2  & 0 \\
0               & a^2+b^2 
\end{array}
\right)=
\left(
\begin{array}{cc}
1  & 0 \\
0  & 1
\end{array}
\right) \, .
\end{equation}
Therefore
\begin{equation}
a^2+b^2 = 1 \, .
\end{equation}
Note that
\begin{equation}
{\rm det}\left(
\begin{array}{cc}
a  & b \\
-b & a
\end{array}
\right) = 1 \, .
\end{equation}
Orthogonality of $ A $ also implies that  blocks of size $ 1 \times 1 $ can be either 1 or~${-1}$. The determinant of the entire matrix $ A $ is the product of determinants of the blocks; the orthogonal transformation of the determinants does not change. Since the determinant of $ Q $ equals one, it follows that the number of $ 1 \times 1 $ blocks containing~$-1$ is even, so that they can be pairwise combined into blocks of the form~(\ref{eq:block}) with $ a = -1 $ , $ b = 0 $. Thus, we can assume that there are blocks of the form (\ref{eq:block}) and blocks with one element~1.

Now we show that there is an exponential representation of the matrix $ A $. To this end, it suffices to demonstrate that such a representation exists in the invariant subspaces, i.e. subspaces of blocks~1 and blocks of the form~(\ref{eq:block}). The first case is trivial, so just consider the blocks~(\ref{eq:block}). We denote such a block as $ B $ and since $ a^2 + b^2 = 1 $, we can write:
\begin{equation}
\label{eq:blockB}
B=\left(
\begin{array}{cc}
\cos\phi  & \sin\phi \\
-\sin\phi & \cos\phi
\end{array}
\right) \, .
\end{equation}
Consider the exponential function $e^{\phi G}$, where the matrix $G$ is of the form:
\begin{equation}
\label{eq:blockG}
G=\left(
\begin{array}{cc}
0   & 1 \\
-1  & 0
\end{array}
\right) \, .
\end{equation}
We expand the exponential function in a power series and note that for the even powers it holds that
\begin{equation}
G^{2n} =(-1)^n E \, ,
\end{equation}
where $E$ is the identity matrix. For the odd powers we have
\begin{equation}
G^{2n+1} =(-1)^n G \, .
\end{equation}
As a result, the power series for the exponential function splits into two series, one of them gives the cosine (even degrees), the other one (odd degrees) gives the sine. So,
\begin{equation}
e^{\phi G}=\cos\phi E + \sin\phi G = B \, .
\end{equation}
Thus, existence of an exponential representation of $ B $, and therefore for $ A $, is proved. 
\begin{equation}
A=e^C = \sum_{n=1}^{\infty} \frac{C^n}{n!}\, ,
\end{equation}
where $ C $ is a matrix of a block-diagonal form, with blocks corresponding to the block given above. Let us multiply the equality by $ W^T $ on the left side and by $ W $ on the right side. Between the factors in powers, we insert identities of the form $ WW^T $. Then
\begin{equation}
Q=W^TAW=\sum_{n=1}^{\infty} \frac{(W^TCW)^n}{n!} = e^{W^TCW}\, .
\end{equation}
The second part of the proposition (that the exponential function of any skew-symmetric matrix is an orthogonal matrix) is based on Proposition~\ref{prop:commute} and the fact that matrix transposition can be transferred to the exponent (the latter may be demonstrated considering the exponential series).

\qed

%
\end{document}